\documentclass{amsart}
\usepackage[T1]{fontenc}
\usepackage[latin1]{inputenc}
\usepackage{fancyhdr}
\usepackage{graphics}
\makeindex

\makeatletter

\providecommand{\LyX}{L\kern-.1667em\lower.25em\hbox{Y}\kern-.125emX\@}

 \theoremstyle{plain}    
 \newtheorem{thm}{Theorem}[section]
 \numberwithin{equation}{section} 
 \numberwithin{figure}{section} 
 \theoremstyle{plain}    
 \newtheorem{cor}[thm]{Corollary} 
 \theoremstyle{plain}    
 \theoremstyle{plain}    
 \newtheorem{prop}[thm]{Proposition} 
 \theoremstyle{definition}
 \newtheorem{defn}[thm]{Definition}
 \theoremstyle{definition}
  
 \theoremstyle{remark}
 \newtheorem{rem}[thm]{Remark}
 \theoremstyle{remark}    
 
 \theoremstyle{remark}    
 \newtheorem*{claim*}{Claim}

\usepackage{verbatim} 
\usepackage{varioref}
\usepackage{amssymb}
\usepackage[mathscr]{eucal}
\usepackage{amsthm}
\usepackage[all]{xy}
\usepackage{hyperref}

\mathcode`\:="603A 

\DeclareSymbolFont{rsfs}{U}{rsfs}{m}{n}
\DeclareSymbolFontAlphabet{\mathrf}{rsfs}

\newcommand{\tmap}[1]{{\mathrm T}_{#1}}

\newcommand{\rank}{\operatorname{rank}}

\newcommand{\glist}[3]{#1_{#2},\dots,#1_{#3}}
\newcommand{\blist}[2]{\glist{#1}{1}{#2}}

\newcommand{\tlist}[2]{\tmap{\blist{#1}{#2}}}

\newcommand{\nth}[1]{$#1^{\mathrm{th}}$}
\newcommand{\tunder}[2]{\tmap{\underbrace{\scriptstyle#1}_{{\text{#2}}}}}

\newcommand{\tunderi}[2]{\tunder{1,\dots,#1,\dots,1}{\nth #2\ position}}

\DeclareFontFamily{U}{manual}{}
\DeclareFontShape{U}{manual}{m}{n}{ <->  manfnt }{}
\newcommand{\manfntsymbol}[1]{%
    {\fontencoding{U}\fontfamily{manual}\selectfont\symbol{#1}}}

\newcommand{\dbend}{\manfntsymbol{127}}
\newcommand{\lhdbend}{\manfntsymbol{126}}
\newcommand{\reversedvideodbend}{\manfntsymbol{0}}
\newcommand{\textdbend}{\text@dbend{\dbend}}
\newcommand{\textlhdbend}{\text@dbend{\lhdbend}}
\newcommand{\textreversedvideodbend}{\text@dbend{\reversedvideodbend}}
\newlength{\dbend@height}
\newcommand{\text@dbend}[1]{%
  \settoheight{\dbend@height}{#1}%
  \divide\dbend@height by 15%
  \multiply\dbend@height by 22%
  \raisebox{\dbend@height}{#1}}

\CompileMatrices

\SelectTips{cm}{}
\newdir{ >}{{}*!/-5pt/@{>}}
\makeatother

\begin{document}


\title{Operads and algebraic homotopy II:\\
Suspensions}

\author{Justin R. Smith}

\subjclass{55R91; Secondary: 18G30}

\keywords{stable homotopy type, minimal models, operads}

\email{jsmith@mcs.drexel.edu}

\urladdr{http://vorpal.mcs.drexel.edu}

\maketitle

\newcommand{\integers}{\mathbb{Z}}

\newcommand{\betabar}{\bar{\beta }}
 
\newcommand{\desusp}{\downarrow }

\newcommand{\susp}{\uparrow }

\newcommand{\cobar}{\mathcal{F}}

\newcommand{\mfrac}{\mathfrak{M}}

\newcommand{\coend}{\mathrm{CoEnd}}

\newcommand{\ainfty}{A_{\infty }}

\newcommand{\coassoc}{\mathrm{Coassoc}}

\newcommand{\trm}{\mathrm{T}}

\newcommand{\tfr}{\mathfrak{T}}

\newcommand{\tabbr}{\hat{\trm }}

\newcommand{\Tabbr}{\hat{\tfr }}

\newcommand{\afr}{\mathfrak{A}}

\newcommand{\homz}{\mathrm{Hom}_{\integers }}

\newcommand{\zend}{\mathrm{End}}

\newcommand{\rs}[1]{\mathrm{R}S_{#1 }}

\newcommand{\highprod}[1]{\bar{\mu }_{#1 }}

\newcommand{\barcs}{\bar{\mathcal{B}}}

\newcommand{\ubarcs}{\mathcal{B}}

\newcommand{\zs}[1]{\mathbb{Z}S_{#1 }}

\newcommand{\homzs}[1]{\mathrm{Hom}_{\integers S_{#1 }}}

\newcommand{\zpi}{\mathbb{Z}\pi }

\newcommand{\D}{\mathfrak{D}}

\newcommand{\cbar}{\hat{\bar{C}}}

\newcommand{\cf}[1]{\mathscr{C}(#1 )}

\newcommand{\cfp}[1]{\mathscr{C}(#1 )^{+}}

\newcommand{\ddelta}{\dot{\Delta }}

\section{Introduction}

The present paper forms a continuation of \cite{Smith:2000}. Since that paper
has not appeared, we begin by summarizing some of its relevant results.

\newdir{ >}{{}*!/-5pt/@{>}}

The easiest way to convey the flavor of this paper's results is with a simple
example. Suppose \( X \) and \( Y \) are pointed, simply-connected, 2-reduced
simplicial sets. There are many topological invariants associated with the chain-complexes
of \( X \) and \( Y \) --- including the coproduct and \emph{\( S_{n} \)-}equivariant
\emph{higher coproducts} (used to define Steenrod operations):
\begin{eqnarray*}
\rs{{n}}\otimes C(X) & \rightarrow  & C(X)^{n}\\
\rs{{n}}\otimes C(Y) & \rightarrow  & C(Y)^{n}
\end{eqnarray*}
for all \( n>1 \), where:

\begin{enumerate}
\item \( \rs{{n}} \) is the bar-resolution of \( \integers  \) over \( \zs{{n}} \).
\item \( (*)^{n} \) denotes the \( n \)-fold tensor product over \( \integers  \)
(with \( S_{n} \) acting by permuting factors).
\end{enumerate}
Now suppose we know (from a purely \emph{algebraic} analysis of these \emph{chain-complexes})
that there exists a chain-map inducing homology isomorphisms
\[
f:C(X)\rightarrow C(Y)\]
and making \begin{equation}\xymatrix{{\rs{n}\otimes C(X)} \ar[r]\ar[d]_{1\otimes f}& {C(X)^n} \ar[d]^{f^n}\\ {\rs{n}\otimes C(Y)} \ar[r]& {C(Y)^n}}\label{dia:intro1}\end{equation}
commute for all \( n>1 \) (requiring exact commutativity is unnecessarily restrictive,
but we assume it to simplify this discussion).

Then our theory asserts that \( X \) and \( Y \) are homotopy equivalent via
a geometric map inducing a map chain-homotopic to \( f \). In fact, our theory
goes further than this: if \( f \) is \emph{any} chain-map making \ref{dia:intro1}
commute for all \( n>1 \), there exists a map of topological realizations
\[
F:|X|\rightarrow |Y|\]
inducing a map equivalent to \( f \) (see \cite[theorems 9.5 and 9.7]{Smith:2000}
for the exact statements). Essentially, Steenrod operations on the \emph{chain-level}
determine: 

\begin{itemize}
\item the homotopy types of spaces and maps,
\item \emph{all} obstructions to topologically realizing chain-maps.
\end{itemize}
I view this work as a generalization of Quillen's characterization of rational
homotopy theory in \cite{Quillen:1969} --- he showed that rational homotopy
types are determined by a commutative coalgebra structures on chain-complexes;
we show that integral homotopy types are determined by coalgebra structures
augmented by higher diagonals.

In \cite[\S~2]{Smith:2000}, we define a category of m-coalgebras, \( \mathfrak{L}_{0} \)
and a localization of it, \( \mathfrak{L} \) by a class of maps called elementary
equivalences. This gives rise to a \emph{homotopy theory} of m-coalgebras in
terms of which we can state this paper's main result:\par{} \vspace{0.3cm}

\emph{Theorem~} \cite{Smith:2000}\emph{:The functor 
\[
\mathscr{C}(*):\underline{\mathrm{Homotop}}_{0}\rightarrow \mathfrak{L}^{+}\]
 defines an equivalence of categories and homotopy theories (in the sense of
\cite{Quillen:1967}), where \( \underline{\mathrm{Homotop}}_{0} \) is the
homotopy category of pointed, simply-connected CW-complexes and continuous maps
and \( \mathfrak{L}^{+}\subset \mathfrak{L} \) is the subcategory of topologically
realizable m-coalgebras. In addition, there exists an equivalence of categories
and homotopy theories
\[
\cf{{*}}:\mathrf {F}\rightarrow \mathfrak{F}^{+}\]
where \( \mathrf {F} \) is the homotopy category of finite, pointed, simply-connected
simplicial sets and \( \mathfrak{F}^{+} \) is the homotopy category of finitely
generated, topologically realizable m-coalgebras in \( \mathfrak{L}_{0} \)
localized with respect to finitely generated equivalences in \( \mathfrak{L}_{0} \).} \par{} \vspace{0.3cm}

In the spirit of our initial statement, one corollary to this result is\par{} \vspace{0.3cm}

\emph{Corollary~9.20 of~\cite{Smith:2000}: Let \( X \) and \( Y \) be pointed,
simply-connected semisimplicial sets and let 
\[
f:\cf{{X}}\rightarrow \cf{{Y}}\]
be a chain-map between canonical chain-complexes. Then \( f \) is topologically
realizable if and only if there exists an m-coalgebra \( C \) over \( \mathfrak{S} \)
and a factorization \( f=f_{\beta }\circ f_{\alpha } \) \[\xymatrix{{\cf{X_1}}\ar[r]^-{f_{\alpha}}&{C}\ar@<.5ex>@{ .>}[r]^-{f_{\beta}}&{\cf{X_2}}\ar@<.5ex>@{ >-}[l]^-{\iota}}\]
where \( f_{\alpha } \) is a morphism of m-coalgebras, \( \iota  \) is an}
elementary equivalence \emph{--- an injection of m-coalgebras with acyclic,
\( \integers  \)-free cokernel --- and \( f_{\beta } \) is a chain map that
is a left inverse to \( \iota  \). If \( X \) and \( Y \) are finite, we
may require \( C \) to be finitely generated.}

\par{} \vspace{0.3cm}The present paper will study \emph{stable homotopy} theory
using similar techniques.

In \S~\ref{sec:suspension}, we define the important concept of suspension of
operads and characterize suspensions of spaces in terms of m-coalgebras. This
leads to the definition of the operad \( \mathfrak{S}_{-\infty } \), which
can be thought of as an infinite desuspension of the operad \( \mathfrak{S} \).
The main result of this section is: \par{} \vspace{0.3cm}

\emph{Corollary~\ref{cor:arbitrarysuspensionspaces}: Let \( X \) and \( Y \)
be pointed, simply-connected semisimplicial sets. Then there exists an integer
\( k\geq 0 \) and a pointed map of topological realizations \( f:|S^{k}X|\rightarrow |S^{k}Y| \)
if and only if there exists a morphism \( g:\mathfrak{p}_{0}^{*}\cfp{{X}}\rightarrow \mathfrak{p}_{0}^{*}\cfp{{Y}} \)
in \( \Sigma ^{-\infty }\mathfrak{L} \) (the category of coalgebras over \( \mathfrak{S}_{-\infty } \)).
In this case, the following diagram commutes in \( \mathfrak{L} \)} \begin{equation}\xymatrix{{\Sigma^k(\mathfrak{U}^k)^*\cfp{X}}\ar@{=}[d] \ar[r]^{\Sigma^kg}& {\Sigma^k(\mathfrak{U}^k)^*\cfp{Y}} \ar@{=}[d]\\ {\cfp{S^kX}} \ar[d]_{\iota_{S^kX}}& {\cfp{S^kY}}\ar[d]^{\iota_{S^kY}} \\ {\cf{\ddelta(|S^kX|)}}\ar[r]_{\cf{\ddelta(f)}} & {\cf{\ddelta(|S^kY|)}}}\end{equation} \emph{where}

\begin{enumerate}
\item \emph{\( \Sigma ^{k} \) denotes suspension of chain-complexes and base-operads.
This is suspension in the sense of dimension-shifting.}
\item \emph{\( \ddelta (*) \) denotes the 2-reduced singular complex functor.}
\item \emph{\( |*| \) denotes topological realization functor.}
\item \emph{\( \iota _{X}:\cf{{X}}\rightarrow \cf{{\ddelta (|X|)}} \) is the canonical
(injective) equivalence.}
\end{enumerate}
\emph{In particular \( X \) and \( Y \) are stably homotopy equivalent if
and only if \( \mathfrak{p}_{0}^{*}\cfp{{X}} \) and \( \mathfrak{p}_{0}^{*}\cfp{{Y}} \)
are equivalent in the category \( \Sigma ^{-\infty }\mathfrak{L} \).}

In future papers, we will apply these results to Boardman's \emph{stable homotopy
category.} The stable homotopy category contains much more than suspensions
of spaces --- for instance, it also includes arbitrary ``formal desuspensions''
of spaces and ``idealized'' space-objects that characterize generalized homology
theories. Even \emph{spaces} behave differently in this category --- as though
already ``infinitely suspended.'' The objects of this category are \emph{spectra}
--- sequences of spaces with maps from the suspension of each term of the sequence
to the next term.

\section{Definitions}

In this section, we recall some of the relevant definitions from \cite{Smith:2000}.

\begin{defn}
\label{def:degmap} Let \( C \) and \( D \) be two graded \( \integers  \)-modules.
A map of graded modules \( f:C_{i}\rightarrow D_{i+k} \) will be said to be
of degree \( k \). 
\end{defn}
\begin{rem}
For instance the \textit{differential} of a chain-complex will be regarded as
a degree \( -1 \) map. 
\end{rem}
We will make extensive use of the Koszul Convention (see~\cite{Gugenheim:1960})
regarding signs in homological calculations:

\begin{defn}
\label{def:basic} \label{def:koszul} If \( f:C_{1}\rightarrow D_{1} \), \( g:C_{2}\rightarrow D_{2} \)
are maps, and \( a\otimes b\in C_{1}\otimes C_{2} \) (where \( a \) is a homogeneous
element), then \( (f\otimes g)(a\otimes b) \) is defined to be \( (-1)^{\deg (g)\cdot \deg (a)}f(a)\otimes g(b) \). 
\end{defn}
\begin{rem}
This convention simplifies many of the common expressions that occur in homological
algebra --- in particular it eliminates complicated signs that occur in these
expressions. For instance the differential, \( \partial _{\otimes } \), of
the tensor product \( \partial _{C}\otimes 1+1\otimes \partial _{D} \).

Throughout this entire paper we will follow the convention that group-elements
act on the left. Multiplication of elements of symmetric groups will be carried
out accordingly --- i.e. 
\[
\left( \begin{array}{cccc}
1 & 2 & 3 & 4\\
2 & 3 & 1 & 4
\end{array}\right) \cdot \left( \begin{array}{cccc}
1 & 2 & 3 & 4\\
4 & 3 & 2 & 1
\end{array}\right) =\]
 result of applying \( \left( \begin{array}{cccc}
1 & 2 & 3 & 4\\
2 & 3 & 1 & 4
\end{array}\right)  \) after applying \( \left( \begin{array}{cccc}
1 & 2 & 3 & 4\\
4 & 3 & 2 & 1
\end{array}\right)  \) or n\( \left( \begin{array}{cccc}
1 & 2 & 3 & 4\\
4 & 1 & 3 & 2
\end{array}\right)  \).

Let \( f_{i} \), \( g_{i} \) be maps. It isn't hard to verify that the Koszul
convention implies that \( (f_{1}\otimes g_{1})\circ (f_{2}\otimes g_{2})=(-1)^{\deg (f_{2})\cdot \deg (g_{1})}(f_{1}\circ f_{2}\otimes g_{1}\circ g_{2}) \).

We will also follow the convention that, if \( f \) is a map between chain-complexes,
\( \partial f=\partial \circ f-(-1)^{\deg (f)}f\circ \partial  \). The compositions
of a map with boundary operations will be denoted by \( \partial \circ f \)
and \( f\circ \partial  \) --- see \cite{Gugenheim:1960}. This convention
clearly implies that \( \partial (f\circ g)=(\partial f)\circ g+(-1)^{\deg (f)}f\circ (\partial g) \).
We will call any map \( f \) with \( \partial f=0 \) a chain-map. We will
also follow the convention that if \( C \) is a chain-complex and \( \susp :C\rightarrow \Sigma C \)
and \( \desusp :C\rightarrow \Sigma ^{-1}C \) are, respectively, the suspension
and desuspension maps, then \( \susp  \) and \( \desusp  \) are both chain-maps.
This implies that the boundary of \( \Sigma C \) is \( -\susp \circ \partial _{C}\circ \desusp  \)
and the boundary of \( \Sigma ^{-1}C \) is \( -\desusp \circ \partial _{C}\circ \susp  \).
\end{rem}
\begin{defn}
\label{r:koszul.5} We will use the symbol \( T \) to denote h transposition
operator for tensor products of chain-complexes \( T:C\otimes D\rightarrow D\otimes C \),
where \( T(c\otimes d)=(-1)^{\dim (c)\cdot \dim (d)}d\otimes c \).
\end{defn}
\begin{prop}
\label{prop:suspisos}Let \( C \) and \( D \) be chain-complexes. Then there
exist isomorphisms
\[
L_{k}=\desusp ^{k}_{C\otimes D}\circ (\susp ^{k}_{C}\otimes 1_{D}):\Sigma ^{-k}C\otimes D\rightarrow \Sigma ^{-k}(C\otimes D)\]
sending \( c\otimes d\in \Sigma ^{-k}C_{i}\otimes D_{j} \) to \( c\otimes d\in \Sigma ^{-k}(C\otimes D)_{i+j} \),
r , \( d\in D_{j} \), and 
\[
M_{k}=\desusp ^{k}_{C\otimes D}\circ (1_{C}\otimes \susp ^{k}_{D}):C\otimes \Sigma ^{-k}D\rightarrow \Sigma ^{-k}(C\otimes D)\]
sending \( c\otimes d\in C_{i}\otimes \Sigma ^{-k}D_{j} \) to \( (-1)^{ik}c\otimes d\in \Sigma ^{-k}(C\otimes D)_{i+j} \),
for \( c\in C_{i} \) and \( d\in \Sigma ^{-k}D_{j}=D_{j+k} \).
\end{prop}
\begin{defn}
\label{def:tmap} Let \( \alpha _{i} \), \( i=1,\dots ,n \) be a sequence
of nonnegative integers whose sum is \( |\alpha | \). Define a set-mapping
 symmetric groups 
\[
\tlist {\alpha }{n}:S_{n}\rightarrow S_{|\alpha |}\]
 as follows:
\begin{enumerate}
\item for \( i \) between 1 and \( n \), let \( L_{i} \) denote the length-\( \alpha _{i} \)
integer sequence: 
\item ,where \( A_{i}=\sum _{j=1}^{i-1}\alpha _{j} \) --- so, for instance, the concatenation
of all of the \( L_{i} \) is the sequence of integers from 1 to \( |\alpha | \); 
\item \( \tlist {\alpha }{n}(\sigma ) \) is the permutation on the integers \( 1,\dots ,|\alpha | \)
that permutes the blocks \( \{L_{i}\} \) via \( \sigma  \). In other words,
 \( \sigma  \) s the permutation 
\[
\left( \begin{array}{ccc}
1 & \dots  & n\\
\sigma (1) & \dots  & \sigma (n)
\end{array}\right) \]
 then \( \tlist {\alpha }{n}(\sigma ) \) is the permutation defined by writing
\[
\left( \begin{array}{ccc}
L_{1} & \dots  & L_{n}\\
L_{\sigma (1)} & \dots  & L_{\sigma (n)}
\end{array}\right) \]
 and regarding the upper and lower rows as sequences length \( |\alpha | \). 
\end{enumerate}
\end{defn}
\begin{rem}
Do not confuse the \( T \)-maps defined here with the transposition map for
tensor products of chain-complexes. We will use the special notation \( T_{i} \)
to represent \( T_{1,\dots ,2,\dots ,1} \), where the 2 occurs in the \( i^{\mathrm{th}} \)
position. The two notations don't conflict since the old notation is never used
in the case when \( n=1 \). Here is an example of the computation of \( \tmap {2,1,3}((1,3,2))=\tmap {2,1,3}\left( \begin{array}{ccc}
1 & 2 & 3\\
3 & 1 & 2
\end{array}\right)  \):\( L_{1}=\{1\}2 \), \( L_{2}=\{3\} \), \( L_{3}=\{4,5,6\} \). The permutation
maps the ordered set \( \{1,2,3\} \) to \( \{3,1,2\} \), so we carry out the
corresponding mapping of the sequences \( \{L_{1},L_{2},L_{3}\} \) to get \( \left( \begin{array}{ccc}
L_{1} & L_{2} & L_{3}\\
L_{3} & L_{1} & L_{2}
\end{array}\right) =\left( \begin{array}{ccc}
\{1,2\} & \{3\} & \{4,5,6\}\\
\{4,5,6\} & \{1,2\} & \{3\}
\end{array}\right) =\left( \begin{array}{cccccc}
1 & 2 & 3 & 4 & 5 & 6\\
4 & 5 & 6 & 1 & 2 & 3
\end{array}\right)  \) (or \( ((1,4)(2,5)(3,6)) \), in cycle notation).
\end{rem}
\begin{defn}
\label{def:operad} A sequence of differential graded \( \mathbb{Z} \)-free
modules, \( \{\mathscr{U}_{i}\} \), will be said to form an \emph{operad} if
they satisfy the following conditions:
\begin{enumerate}
\item there exists a \emph{unit map} (defined by the commutative diagrams below) 
\[
\eta :\mathbb{Z}\rightarrow \mathscr{U}_{1}\]

\item for all \( i>1 \), \( \mathscr{U}_{i} \) is equipped with a left action of
\( S_{i} \), the symmetric group. 
\item for all \( k\geq 1 \), and \( i_{s}\geq 0 \) there are maps 
\[
\gamma :\mathscr{U}_{i_{1}}\otimes \cdots \otimes \mathscr{U}_{i_{k}}\otimes \mathscr{U}_{k}\rightarrow \mathscr{U}_{i}\]
 where \( i=\sum _{j=1}^{k}i_{j} \).

The \( \gamma  \)-maps must satisfy the following conditions:

\end{enumerate}
\end{defn}
\begin{description}
\item [Associativity]the following diagrams commute, where \( \sum j_{t}=j \), \( \sum i_{s}=i \),
and \( g_{\alpha }=\sum _{\ell =1}^{\alpha }j_{\ell } \) and \( h_{s}=\sum _{\beta =g_{s-1}+1}^{g_{s}}i_{\beta } \):
 \[\xymatrix@C+20pt{{\left(\bigotimes_{s=1}^{j}\mathscr{U}_{i_{s}}\right)\otimes\left(\bigotimes_{t=1}^{k}\mathscr{U}_{j_{t}}\right) \otimes\mathscr{U}_{k}}\ar[r]^-{\text{Id}\otimes\gamma}\ar[dd]_{\text{shuffle}}&{\left(\bigotimes_{s=1}^{j}\mathscr{U}_{i_{s}}\right)\otimes\mathscr{U}_{j}}\ar[d]^{\gamma}\\
&{\mathscr{U}_{i}}\\{\left(\left(\bigotimes_{q=1}^{j_{t}}\mathscr{U}_{i_{g_{t-1}+q}}\right)\otimes\bigotimes_{t=1}^{k}\mathscr{U}_{j_{t}}\right) \otimes\mathscr{U}_{k}}\ar[r]_-{(\otimes_{t}\gamma)\otimes\text{Id}}&{\left(\bigotimes_{t=1}^{k}\mathscr{U}_{h_{k}}\right) \otimes\mathscr{U}_{k}}\ar[u]_{\gamma}}\]

\item [Units]the following diagrams commute:  \[\begin{array}{cc}\xymatrix{{{\integers}^{k}\otimes\mathscr{U}_{k}}\ar[r]^{\cong}\ar[d]_{{\eta}^{k}\otimes\text{Id}}&{\mathscr{U}_{k}}\\
{{\mathscr{U}_{1}}^{k}\otimes{\mathscr{U}_{k}}}\ar[ur]_{\gamma}&}&\xymatrix{{\mathscr{U}_{k}\otimes\integers}\ar[r]^{\cong}\ar[d]_{\text{Id}\otimes\eta}&{\mathscr{U}_{k}}\\
{\mathscr{U}_{k}\otimes\mathscr{U}_{1}}\ar[ur]_{\gamma}&}\end{array}\] 
\item [Equivariance]the following diagrams commute: \[\xymatrix@C+20pt{{\mathscr{U}_{j_{1}}\otimes\cdots\otimes\mathscr{U}_{j_{k}}\otimes\mathscr{U}_{k}}\ar[r]^-{\gamma}\ar[d]_{\sigma^{-1}\otimes\sigma}&{\mathscr{U}_{j}}\ar[d]^{\tmap{j_{1},\dots,j_{k}}(\sigma)}\\
{\mathscr{U}_{j_{\sigma(1)}}\otimes\cdots\otimes\mathscr{U}_{j_{\sigma(k)}}\otimes\mathscr{U}_{k}}\ar[r]_-{\gamma}&{\mathscr{U}_{j}}}\]
where \( \sigma \in S_{k} \), and the \( \sigma ^{-1} \) on the left permutes
the factors \( \{\mathscr{U}_{j_{i}}\} \) and the \( \sigma  \) on the right
simply acts on \( \mathscr{U}_{k} \). See \ref{def:tmap} for a definition
of \( \tmap {j_{1},\dots ,j_{k}}(\sigma ) \). \[\xymatrix@C+20pt{{\mathscr{U}_{j_{1}}\otimes\cdots\otimes\mathscr{U}_{j_{k}}\otimes\mathscr{U}_{k}}\ar[r]^-{\gamma}\ar[d]_{\tau_{1}\otimes\cdots\tau_{k}\otimes\text{Id}}&{\mathscr{U}_{j}}\ar[d]^-{\tau_{1}\oplus\cdots\oplus\tau_{k}}\\
{\mathscr{U}_{j_{\sigma(1)}}\otimes\cdots\otimes\mathscr{U}_{j_{\sigma(k)}}\otimes\mathscr{U}_{k}}\ar[r]_-{\gamma}&{\mathscr{U}_{j}}}\] where
\( \tau _{s}\in S_{j_{s}} \) and \( \tau _{1}\oplus \cdots \oplus \tau _{k}\in S_{j} \)
is the block sum.
\end{description}
\begin{rem}
The alert reader will notice a discrepancy between our definition of operad
and that in \cite{Kriz-May} (on which it was based). The difference is due
to our using operads as parameters for systems of \emph{maps}, rather than \( n \)-ary
operations. We, consequently, compose elements of an operad as one composes
\emph{maps}, i.e. the second operand is to the \emph{left} of the first. This
is also why the symmetric groups act on the \emph{left} rather than on the right. 
\end{rem}
We will frequently want to think of operads in other terms:

\begin{defn}
\label{def:operadcomps} Let \( \mathscr{U} \) be an operad as defined above.
Given \( k_{1}\geq k_{2}>0 \), define the \( i^{\mathrm{th}} \) \emph{composition}

\[
\circ _{i}:\mathscr{U}_{k_{2}}\otimes \mathscr{U}_{k_{1}}\rightarrow \mathscr{U}_{k_{1}+k_{2}}\]
 as the composite \begin{multline}\underbrace{\integers\otimes\cdots\otimes\integers\otimes\mathscr{U}_{k_{2}}\otimes\integers\otimes\cdots\otimes\integers}_{\text{\(i^{\text{th}}\)factor}}\otimes\mathscr{U}_{k_{1}}\\ \to\underbrace{\mathscr{U}_{1}\otimes\cdots\otimes\mathscr{U}_{1}\otimes\mathscr{U}_{k_{2}}\otimes\mathscr{U}_{1}\otimes\cdots\otimes\mathscr{U}_{1}}_{\text{\(i^{\text{th}}\)factor}}\otimes\mathscr{U}_{k_{1}}\to\mathscr{U}_{k_{1}+k_{2}-1}\end{multline}
where the final map on the right is \( \gamma  \). 

These compositions satisfy the following conditions, for all \( a\in \mathscr{U}_{n} \),
\( b\in \mathscr{U}_{m} \), and \( c\in \mathscr{U}_{t} \):
\begin{description}
\item [Associativity]\( (a\circ _{i}b)\circ _{j}c=a\circ _{i+j-1}(b\circ _{j}c) \)
\item [Commutativity]\( a\circ _{i+m-1}(b\circ _{j}c)=(-1)^{mn}b\circ _{j}(a\circ _{i}c) \)
\item [Equivariance]\( a\circ _{\sigma (i)}(\sigma \cdot b)=\tunderi {n}{i}(\sigma )\cdot (a\circ _{i}b) \)
\end{description}
\end{defn}
\begin{rem}
In \cite{Smith:1994}, I originally \emph{defined} operads (or formal coalgebras)
in terms of these compositions. It turned out that I'd recapitulated the historical
sequence of events: operads were originally defined this way and called \emph{composition
algebras.} I am indebted to Jim Stasheff for pointing this out to me. 

Given this definition of operad, we recover the \( \gamma  \) map in \ref{def:operad}
by setting: 
\[
\gamma (u_{i_{1}}\otimes \cdots \otimes u_{i_{k}})=u_{i_{1}}\circ _{1}\cdots \circ _{k-1}u_{i_{k}}\circ _{k}u_{k}\]
 (where the implied parentheses associate to the right). It is left to the reader
to verify that the two definitions are equivalent (the commutativity condition,
here, is a special case of the equivariance condition).
\end{rem}
Morphisms of operads are defined in the obvious way:

\begin{defn}
\label{def:operadmorphism} Given two operads \( \mathscr{U} \) and \( \mathscr{V} \),
a \emph{morphism} 
\[
f:\mathscr{U}\rightarrow \mathscr{V}\]
 is a sequence of chain-maps 
\[
f_{i}:\mathscr{U}_{i}\rightarrow \mathscr{V}_{i}\]
 commuting with all the diagrams in \ref{def:operad} or (equivalently) preserving
the composition operations in \ref{def:operadcomps}.
\end{defn}
Now we give some examples:

\begin{defn}
\label{def:mathfrakS0}The operad \( \mathfrak{S}_{0} \) is defined via
\end{defn}
\begin{enumerate}
\item Its \( n^{\mathrm{th}} \) component is \( \zs{{n}} \) --- a chain-complex
concentrated in dimension \( 0 \). 
\item The composition operations are defined by
\[
\gamma (\sigma _{i_{1}}\otimes \cdots \otimes \sigma _{i_{k}}\otimes \sigma )=\sigma _{i_{\sigma (k)}}\oplus \cdots \oplus \sigma _{i_{\sigma (k)}}\circ \tmap {i_{1},\dots ,i_{k}}(\sigma )\]

\end{enumerate}
\begin{rem}
This was denoted \( \mathrf {M} \) in \cite{Kriz-May}.
\end{rem}
Verification that this satisfies the required identities is left to the reader
as an exercise.

\begin{defn}
\label{def:sfrakfirstmention}Let \( \mathfrak{S} \) denote the \( E_{\infty } \)-operad
with components
\[
\rs{{n}}\]
--- the bar resolutions of \( \integers  \) over \( \zs{{n}} \) for all \( n>0 \).
This is an important operad and \cite[\S~2.3]{Smith:2000} is devoted to it.
\end{defn}
\begin{rem}
This is the result of applying the ``unreduced bar construction'' to the previous
example.
\end{rem}
\begin{defn}
\label{def:coassoc}\( \coassoc  \) is an operad defined to have one basis
element \( \{b_{i}\} \) for all integers \( i\geq 0 \). Here the rank of \( b_{i} \)
is \( i \) and the degree is 0 and the these elements satisfy the composition-law:
\( b_{i}\circ _{\alpha }b_{j}=b_{i+j-1} \) regardless of the value of \( \alpha  \),
which can run from \( 1 \) to \( j \). The differential of this operad is
identically zero.
\end{defn}
Now we define two important operads associated to any \( \integers  \)-module.

\begin{defn}
\label{def:coend} Let \( C \) be a DGA-module with augmentation \( \epsilon :C\rightarrow \mathbb{Z} \),
and with the property that \( C_{0}=\mathbb{Z} \). Then the \emph{Coendomorphism}
operad, \( \coend (C) \), is defined to be the operad with: 
\begin{enumerate}
\item component of \( \rank i=\homz (C,C^{i}) \), with the differential induced by
that of \( C \) and \( C^{i} \). The dimension of an element of \( \homz (C,C^{i}) \)
(for some \( i \)) is defined to be its degree as a map. 
\item The \( \mathbb{Z} \)-summand is generated by one element, \( e \), of rank
0. 
\end{enumerate}
\end{defn}
\begin{rem}
Both operads are unitary --- their common identity element is the identity map
\( \mathrm{id}\in \homz (C,C) \). One motivation for operads is that they model
the iterated coproducts that occur in \( \coend (*) \). We will use operads
as an algebraic framework for defining other constructs that have topological
applications.
\end{rem}
In like fashion, we define the \textit{endomorphism} operad:

\begin{defn}
\label{def:end} If \( C \) is a DGA-module, the \emph{endomorphism operad},
\( \zend (C) \) is defined to have components 
\[
\homz (C^{n},C)\]
 and compositions that coincide with endomorphism compositions.
\end{defn}
Now we consider a special class of operads that play a crucial role in the sequel:

\begin{defn}
\label{def:einfinity} An operad \( \mathscr{U}=\{\mathscr{U}_{n}\} \) will
be called an \( E_{\infty } \)-operad if \( \mathscr{U}_{n} \) is a \( \zs{n} \)-free
resolution of \( \mathbb{Z} \) for all \( n>0 \). 
\end{defn}
Given these definitions, we can define:

\begin{defn}
\label{def:operadcomodule} Let \( \mathscr{U} \) be an operad and let \( C \)
be a DG-module equipped with a morphism (of operads) 
\[
f:\mathscr{U}\rightarrow \coend (C)\]
 Then \( C \) is called a \emph{coalgebra} over \( \mathscr{U} \) with structure
map \( f \). If \( C \) is equipped with a morphism of operads 
\[
f:\mathscr{U}\rightarrow \zend (C)\]
 then \( C \) is called a \emph{algebra} over \( \mathscr{U} \) with structure
map \( f \). 
\end{defn}
\begin{rem}
A coalgebra, \( C \), over an operad, \( \mathscr{U} \), is a sequence of
maps 
\[
f_{n}:\mathscr{U}\otimes C\rightarrow C^{n}\]
 for all \( n>0 \), where \( f_{n} \) is \( \zs{n} \)-equivariant. These
maps are related in the sense that they fit into commutative diagrams:  \[\xymatrix@d@R+20pt{{\mathscr{U}_{n}\otimes\mathscr{U}_{m}\otimes  C}\ar[r]^-{\circ_{i}}&{\mathscr{U}_{n+m-1}\otimes 
C}\ar[r]^-{f_{n+m-1}}&{C^{n+m-1}}\\ 
{\mathscr{U}_{n}\otimes\mathscr{U}_{m}\otimes C}\ar[r]_-{1\otimes 
f_{m}}\ar@{=}[u]&{\mathscr{U}_{n}\otimes 
C^{m}}\ar[r]_-{V_{i-1}}&{C^{i-1}\otimes\mathscr{U}_{n}\otimes C\otimes 
C^{m-i}}\ar[u]_-{1\otimes\dots\otimes f_{n}\otimes\dots\otimes1}}
\]  for all \( n,m\geq 1 \) and \( 1\leq i\leq m \). Here \( V:\mathscr{U}_{n}\otimes C^{m}\rightarrow C^{i-1}\otimes \mathscr{U}_{n}\otimes C\otimes C^{m-i} \)
is the map that shuffles the factor \( \mathscr{U}_{n} \) to the right of \( i-1 \)
factors of \( C \). In other words: The abstract composition-operations in
\( \mathscr{U} \) exactly correspond to compositions of maps in \( \{\homz (C,C^{n})\} \).
We exploit this behavior in applications of coalgebras over operads, using an
explicit knowledge of the algebraic structure of \( \mathscr{U} \).
\end{rem}
In very simple cases, one can explicitly describe the maps defining a coalgebra
over an operad:

\begin{defn}
\label{def:unitinterval}Define the coalgebra, \( I \) --- \emph{the unit interval}
--- over \( \mathfrak{S} \) via:
\end{defn}
\begin{enumerate}
\item Its \( \integers  \)-generators are \( \{p_{0},p_{1}\} \) in dimension \( 0 \)
and \( q \) in dimension \( 1 \), and its adjoint structure map is \( r_{n}:\rs{{n}}\otimes I\rightarrow I^{n} \),
for all \( n>1 \).
\item The coproduct is given by \( r_{2}([\, ]\otimes p_{i})=p_{i}\otimes p_{i} \),
\( i=0,1 \), and \( r_{2}([\, ]\otimes q)=p_{0}\otimes q+q\otimes p_{1} \).
\item The higher coproducts are given by \( r_{2}([(1,2)]\otimes q=q\otimes q \),
\( r_{2}([(1,2)]\otimes p_{i})=0 \), \( i=0,1 \) and \( r_{2}(a\otimes I)=0 \),
where \( a\in \rs{{2}} \) has dimension \( >1 \).
\end{enumerate}
\begin{rem}
This, coupled with the operad-identities in \( \mathfrak{S} \) suffice to define
the coalgebra structure of \( I \) in all dimensions and for all degrees.
\end{rem}
\begin{prop}
Coassociative coalgebras are precisely the coalgebras over \( \coassoc  \). 
\end{prop}
\begin{rem}
There are some subtleties to this definition, however. It is valid if we regard
\( \coassoc  \) as a non-\( \Sigma  \) operad. If we regard it as an operad
with \emph{trivial} symmetric group action, then we have defined coassociative,
\emph{cocommutative} coalgebras. 
\end{rem}
We can define tensor products of operads:

\begin{defn}
Let \( \mathscr{U}_{1} \) and \( \mathscr{U}_{2} \) be operads. Then \( \mathscr{U}_{1}\otimes \mathscr{U}_{2} \)
is defined to have: 
\begin{enumerate}
\item component of \( \rank i=(\mathscr{U}_{1})_{i}\otimes (\mathscr{U}_{2})_{i} \),
where \( (\mathscr{U}_{1})_{i} \) and \( (\mathscr{U}_{2})_{i} \) are, respectively,
the components of \( \rank i \) of \( \mathscr{U}_{1} \) and \( \mathscr{U}_{2} \); 
\item composition operations defined via \( (a\otimes b)\circ _{i}(c\otimes d)=(-1)^{\dim (b)\dim (c)}(a\circ _{i}c)\otimes (b\circ _{i}d) \),
for \( a,c\in \mathscr{U}_{1} \) and \( b,d\in \mathscr{U}_{2} \). 
\end{enumerate}
\end{defn}
We conclude this section with 

\begin{defn}
\label{def:operadcoproduct}An operad, \( \mathfrak{R} \), will be called an
\emph{operad-coalgebra} if there exists a co-associative morphism of operads
\[
\Delta :\mathfrak{R}\rightarrow \mathfrak{R}\otimes \mathfrak{R}\]

\end{defn}
\begin{rem}
Operad-coalgebras are important in certain homotopy-theoretic contexts --- for
instance in the study of the bar and cobar constructions. The \( \mathfrak{S} \)-operad
is a particularly important operad of this type.
\end{rem}

\section{Suspensions\label{sec:suspension}}

\subsection{Chain-complexes and operads}

In this section, we will compute the effect of a suspension in some detail.
We begin with a word on notation. There is the well-known concept of suspension
of a chain-complex: this is just the given chain-complex with dimensions shifted
up or down. In \cite{Smith:2000}, we defined the concept of suspension of an
m-coalgebra that corresponded to the topological notion of suspension. Now the
problem: although m-coalgebras are also chain-complexes, their suspensions as
\emph{chain-complexes} are \emph{not} the same as their suspensions as \emph{m-coalgebras}. 

We will consequently use the notation \( SC \) for the supension of \( C \)
as an \emph{m-coalgebra} (this is the suspension that has topological significance).

\begin{defn}
\label{def:supensionspecialoperad}Let \( C=\Sigma \integers  \) be the chain-complex
concentrated in dimension \( 1 \) and define \( \mathrm{Susp}=\coend (C) \),
the \emph{suspension operad.} In like fashion, we can define \( \mathrm{Susp}^{-1}=\coend (\Sigma ^{-1}\integers ) \),
the \emph{desuspension} operad.
\end{defn}
\begin{rem}
Note that \( \mathrm{Susp}_{n} \) is a chain-complex concentrated in dimension
\( n-1 \) equal to \( \Sigma ^{n-1}\integers  \), with compositions defined
by
\[
s_{n}\circ _{i}s_{m}=(-1)^{(i-1)(n-1)}s_{n+m-1}\]
where \( s_{n}\in \mathrm{Susp}_{n} \) is the canonical generator. The \( S_{n} \)-action
on \( s_{n} \) is defined by
\[
\sigma \cdot s_{n}=(-1)^{\mathrm{parity}(\sigma )}s_{n}\]
where \( \sigma \in S_{n} \).

Note that there are two distinct, but canonically isomorphic, ways of defining
\( \mathrm{Susp} \). We have selected the definition that corresponds to composing
maps \( f,g\in \coend (C) \) to have \( f\circ _{i}g \) mean ``the map \( g \)
\emph{followed by} the map \( f \).''

This definition was motivated by
\end{rem}
\begin{prop}
\label{prop:suspcoend}If \( C \) is a chain-complex
\[
\coend (\Sigma C)=\mathrm{Susp}\otimes \coend (C)\]

\end{prop}
\begin{rem}
Note that the suspensions here are those of \emph{chain-complexes} (i.e., simple
dimension-shifting) rather than the (topologically motivated) notion of \emph{m-coalgebra}
suspension defined in \cite[3.36]{Smith:2000}. We will connect these two notions
later in this section.

For the component of degree \( n \), the isomorphisms in question are precisely
\( \homz (\desusp ,\underbrace{\desusp \otimes \cdots \otimes \desusp }_{n\, \mathrm{times}}) \),
which lower dimension by \( n-1 \).
\end{rem}
\begin{proof}
Note that \( \Sigma C=\Sigma \integers \otimes C \). The conclusion follows
from the injective operad-morphism (see \cite[Proposition 3.4]{Smith:2000})
\[
\mathfrak{E}:\coend (\Sigma \integers )\otimes \coend (C)\rightarrow \coend (\Sigma \integers \otimes C)\]
which is an isomorphism in this case (since \( D=\mathrm{Susp} \) has finitely
generated components). 
\end{proof}
We can use this to define the concept of \emph{suspension of an operad:}

\begin{defn}
\label{def:suspoperad}\label{prop:suspoperadtensorproductsusp}Given an operad
\( A \), we define its suspension \( \Sigma A \) to be the operad \( \mathrm{Susp}\otimes A \)
and its \emph{desuspension,} \( \Sigma ^{-1}A \), to be the operad \( \mathrm{Susp}^{-1}\otimes A \).
The \emph{suspension isomorphism} \( \mathscr{I}:A\rightarrow \Sigma A \) is
a set of isomorphisms 
\[
\mathscr{I}_{n}:A_{n}\rightarrow (\Sigma A)_{n}\]
of degree \( n-1 \) (i.e., \( \mathscr{I}_{n} \) raises dimension by \( n-1) \)
for all \( n>0 \) such that:
\begin{enumerate}
\item The identity
\[
\mathscr{I}_{n}(g\cdot x)=(-1)^{\mathrm{parity}(g)}g\cdot \mathscr{I}_{n}(x)\]
holds for all \( n>0 \), all \( x\in A_{n} \) and all \( g\in S_{n} \).
\item The identity
\[
\mathscr{I}_{n+m-1}(a\circ _{i}b)=(-1)^{(m-1)\dim a+(n-1)(i-1)}\mathscr{I}_{n}(a)\circ _{i}\mathscr{I}_{m}(b)\]
holds for all \( n,m>0 \), all \( a\in A_{n} \) and \( b\in A_{m} \), and
all \( 0\leq i\leq m \).
\end{enumerate}
\end{defn}
Although suspension-isomorphisms are not operad-morphisms (they do not even
preserve dimension) they clearly define a functor from the category of operads
to itself.

\begin{defn}
Let \( k\geq 0 \) be an integer and define define \( \Sigma ^{-k}\mathfrak{L} \)
to be the category of DG-coalgebras over the operad \( \Sigma ^{-k}\mathfrak{S} \)
that are \( k \)-fold desuspensions (in the \( \Sigma  \)-sense) of m-coalgebras.
Define \( \Sigma ^{-k}\mathfrak{L}=\Sigma ^{-k}\mathfrak{L}_{0}[H^{-1}] \)
where \( H \) is the class of morphisms of coalgebras over \( \Sigma ^{-k}\mathfrak{S} \)
inducing \emph{homology isomorphisms} of underlying chain-complexes.
\end{defn}
By \emph{definition}

\begin{prop}
\label{prop:desuspensioninj}Desuspension defines isomorphisms of categories
\[
\desusp ^{k}:\mathfrak{L}\rightarrow \Sigma ^{-k}\mathfrak{L}\]
for all \( k\geq 0 \).
\end{prop}

\subsection{m-coalgebras}

For the remainder of this section, we will focus on m-coalgebras over the operad
\( \mathfrak{S} \) and consider suspensions of m-coalgebras as defined in \cite[3.36]{Smith:2000}. 

There are important differences between suspending m-coalgebras and arbitrary
coalgebras over an operad. For instance m-coalgebras contain a sub-coalgebra
isomorphic to the \emph{trivial m-coalgebra,} \( \integers  \) --- the \emph{basepoint}.
In addition, there exists a canonical surjection of coalgebras
\[
C\rightarrow \integers \]
 that is a left-inverse to inclusion of the basepoint. The basepoint does \emph{not}
participate in suspension. 

\begin{defn}
Let \( k \) be an integer \( \geq 0 \), and let \( C \) be an object of \( \Sigma ^{-k}\mathfrak{L}_{0} \).
Then we define \( C^{+}=C/\Sigma ^{-k}\integers  \), where \( \Sigma ^{-k}\integers =C_{-k} \)
is the (desuspended) basepoint.
\end{defn}
\begin{rem}
By \cite[2.35]{Smith:2000}, this quotient is a coalgebra over \( \Sigma ^{-k}\mathfrak{S} \).
\end{rem}
\begin{defn}
\label{def:reducedmcoalgebra}Let \( k \) be an integer \( \geq 0 \), and
let \( C \) be an object of \( \Sigma ^{-k}\mathfrak{L}_{0} \) with structure-map
\[
a:\Sigma ^{-k}\mathfrak{S}\rightarrow \coend (C)\]
We will call \( C \) \emph{reduced} if its coproduct is trivial -- i.e.
\[
\Delta (x)=a(\Sigma ^{-k}[\, ])(x)=\Sigma ^{-k}1\otimes x+x\otimes \Sigma ^{-k}1\]
for all \( x\in C \), where \( 1\in C_{-k} \) is the (desuspended) basepoint.
\end{defn}
The following is immediate:

\begin{prop}
Let \( k \) be an integer \( \geq 0 \), and let \( C \) and \( D \) be reduced
objects of \( \Sigma ^{-k}\mathfrak{L}_{0} \). Then \( C\cong D \) if and
only if \( C^{+}\cong D^{+} \).
\end{prop}
Now we are in a position to begin to calculate the \( S\wedge * \)-suspension
of m-coalgebras.

\begin{prop}
Let \( [0,1] \) denote the standard 1-simplex, with vertices \( \{[0],[1]\} \),
and let \( \mathrm{Susp} \) be the operad defined in \ref{def:supensionspecialoperad}.
\emph{Then there exists a surjective operad-morphism
\begin{equation}
\label{eq:knmaps}
\tau :\coend (\mathscr{C}([0,1]))\rightarrow \mathrm{Susp}
\end{equation}
}
\end{prop}
\begin{proof}
This is just the morphism of \( \coend (*) \)-operads induced by the chain-map
\[
\mathscr{C}([0,1])\rightarrow \mathscr{C}([0,1])/\mathscr{C}([0,1])_{0}=\Sigma \integers \]
 
\end{proof}
\begin{prop}
\label{prop:vmapsdef}There exists a surjective operad-morphism
\[
\mathfrak{V}:\mathfrak{S}\rightarrow \Sigma \mathfrak{S}\]
 that makes the following diagram commute \[\xymatrix{&{\Sigma \mathfrak{S} }\ar[rd]^-{\mathscr{I} \otimes 1\circ \Delta\circ \mathscr{I}^{-1}}&{}\\ {\mathfrak{S}}\ar[ur]^-{\mathfrak{V}}\ar[rd]_-{\Delta}& {} & {\Sigma\mathfrak{S} \otimes \mathfrak{S}} \\ {} & {\mathfrak{S} \otimes \mathfrak{S}} \ar[ru]_-{\mathfrak{V} \otimes 1}& {}}\]
\end{prop}
\begin{rem}
The morphism \( \mathfrak{V} \) coincides with the \emph{determinant map} defined
in \cite{Smith:2000}. Note that \cite[Proposition 2.26]{Smith:2000} does not
apply since \( \Sigma \mathfrak{S} \) and \( \Sigma ^{-1}\mathfrak{S} \) are
not \( E_{\infty } \)-operads. Although the \( \{\mathfrak{V}_{n}\} \) do
not induce homology isomorphisms, they do define cohomology classes
\[
\{\alpha _{n}\in H^{n-1}(S_{n},\integers )\}\]
for all \( n>0 \).
\end{rem}
\begin{proof}
Identify \( \Sigma \mathfrak{S} \) with the operad \( \mathrm{Susp}\otimes \mathfrak{S} \),
by \ref{prop:suspoperadtensorproductsusp}. Now form the composite \begin{equation} \xymatrix{{\mathfrak{S}}\ar[r]^-{\Delta}&{\mathfrak{S}\otimes \mathfrak{S}}\ar[r]^-{u\otimes 1}& {\coend(U)\otimes \mathfrak{S} }\ar[r]^-{\tau\otimes 1}& {\mathrm{Susp}\otimes\mathfrak{S} = \Sigma\mathfrak{S}} }\label{eq:suspendmathfrakS}\end{equation}
where 
\begin{enumerate}
\item \( \Delta :\mathfrak{S}\rightarrow \mathfrak{S}\otimes \mathfrak{S} \) is the
canonical coproduct defined in \cite{Smith:2000},
\item the morphism \( \tau  \) is defined in \ref{eq:knmaps}, 
\item \( U=\mathscr{C}([0,1]) \), 
\item \( u:\mathfrak{S}\rightarrow \coend (U) \) is the structure map of the m-structure
of \( \mathscr{C}([0,1]) \). 
\end{enumerate}
\end{proof}
We will use this to describe the m-structure of suspensions.

\begin{prop}
\label{prop:geosuspchainsusp}Let \( C \) be an m-coalgebra over \( \mathfrak{S} \).
Then 
\[
\Sigma C^{+}=(SC)^{+}\]
 as chain-complexes.
\end{prop}
\begin{proof}
Recall that
\[
SC=\mathscr{C}(S^{1})\wedge C=\mathscr{C}(I)\otimes C/([0]\otimes C^{+}+[1]\otimes C^{+}\otimes +[0,1]^{+}\otimes p)\]
 where \( p\in C_{0}=\integers  \) is the basepoint and \( [0] \) is the basepoint
of \( [0,1] \). This means that the basepoint of \( SC \) is \( [0]\otimes p \)
and \( (SC)^{+}=[0,1]\otimes C^{+}=\Sigma \integers \otimes C^{+}=\Sigma C^{+} \).
\end{proof}
The following result relates our two distinct notions of suspension

\begin{thm}
\label{th:operadsuspensiondiagram}Let \( C \) be an m-coalgebra over \( \mathfrak{S} \)
with structure map
\[
a_{C}:\mathfrak{S}\rightarrow \coend (C)\]
whose suspension (as an m-coagebra) is \( \Sigma C \) with structure map
\[
a_{SC}:\mathfrak{S}\rightarrow \coend (SC)\]
Then \( \Sigma C^{+}=(SC)^{+} \). Suppose that the corresponding reduced m-coalgebras
have the structure maps
\[
a^{+}_{C}:\mathfrak{S}\rightarrow \coend (C^{+})\]
and
\[
a^{+}_{SC}:\mathfrak{S}\rightarrow \coend ((SC)^{+})=\coend (\Sigma C^{+})\]
respectively. Then the following diagram commutes \begin{equation}\xymatrix{ {\coend((SC)^+)} \ar@{=}[r]& {\coend(\Sigma C^+)} & {\coend(C^+)}\ar[l]_-{\mathscr{I}}\\ {\mathfrak{S}} \ar[u]^-{a^+_{SC}} \ar[r]_{\mathfrak{V}}& {\Sigma \mathfrak{S} }\ar[u]_-{\Sigma a^+_C} & {\mathfrak{S}}\ar[l]^-{\mathscr{I}}\ar[u]_{a_C^+}}\label{eq:bigsuspcommut}\end{equation} where
all vertical maps are operad-morphisms and the the \( \mathscr{I} \)-maps are
suspension-isomorphisms of operads (see \ref{def:suspoperad}).
\end{thm}
\begin{proof}
This follows from \ref{prop:vmapsdef} and the definition of the m-structure
of the (geometric) suspension. Let \( C \) be an m-coalgebra over \( \mathfrak{S} \)
with structure map
\[
a:\mathfrak{S}\rightarrow \coend (C)\]
 and let the structure map of the unit interval \( I=\mathscr{C}([0,1]) \)
be
\[
u:\mathfrak{S}\rightarrow \coend (I)\]
Then the suspension has a structure map that is the composite \begin{equation} \xymatrix{{\mathfrak{S}}\ar@{=}[r]\ar[d]_-{\Delta}&{\mathfrak{S}}\ar[d]^-{\Delta} \\ {\mathfrak{S} \otimes \mathfrak{S}}\ar[d]_-{u\otimes 1}\ar@{=}[r]& {\mathfrak{S} \otimes \mathfrak{S}}\ar[d]^{u \otimes a}\\ {\coend(I) \otimes \mathfrak{S} }\ar[r]^-{1\otimes a }\ar@{=}[d]& { \coend(I)\otimes \coend(C) }\ar[d]^{\mathfrak{E} } \\ {\coend(I)\otimes \mathfrak{S} } \ar[r]^-{\mathfrak{E}\circ 1\otimes a}\ar[d]_-{\{\tau \otimes 1\}} &{\coend(I\otimes C )}\ar[d]^-{q} \\ {\mathrm{Susp}\otimes \mathfrak{S}} \ar@{=}[d]\ar[r]^-{(\mathfrak{E}\circ 1\otimes a)^*} & {\coend((SC)^+)}\ar@{=}[d] \\ {\Sigma \mathfrak{S}} \ar[r]^-{\Sigma a^+} & {\coend(\Sigma C^+)} \\ } \end{equation}Where
\begin{enumerate}
\item \( \mathfrak{E}:\coend (I)\otimes \coend (C)\rightarrow \coend (I\otimes C) \)
is the operad morphism defined in \cite[Proposition 3.4]{Smith:2000}
\item \( q:\coend (I\otimes C)\rightarrow \coend ((SC)^{+}) \) is induced by the
quotient map \( I\otimes C\rightarrow I\otimes C/([0]\otimes C^{+}+[1]\otimes C^{+}+[0,1]^{+}\otimes p+[0]\otimes p)=(SC)^{+} \)
\end{enumerate}
Now observe that the left column is \( \mathfrak{V} \) and that the right is
the structure map of \( (SC)^{+} \).

\end{proof}
\begin{cor}
\label{cor:suspensiondown}Let \( C \) and \( D \) be m-coalgebras over \( \mathfrak{S} \)
with reduced structure morphisms
\[
a_{C}:\mathfrak{S}\rightarrow \coend (C^{+})\]
and
\[
a_{D}:\mathfrak{S}\rightarrow \coend (D^{+})\]
respectively. Then a morphism \( f:SC\rightarrow D \) gives rise to a commutative
diagram \begin{equation}\xymatrix@C+10pt{ {C^+} \ar[r]^-{\Sigma^{-1}f^+}& {\Sigma^{-1}D^+} \\ {\mathfrak{S}} \ar[u]^{a_C}& {\Sigma^{-1}\mathfrak{S}}\ar[u]_{\Sigma^{-1}a_D} \ar[l]^{\mathfrak{U}}}\label{dia:smallsuspensionsquare}\end{equation}
\end{cor}
\begin{proof}
Just desuspend the left square of diagram \ref{eq:bigsuspcommut} (and interchange
its columns) to get \[\xymatrix{ {\coend(C^+)}& {\Sigma^{-1}\coend((SC)^+)} \ar@{=}[l]\\{\mathfrak{S} }\ar[u]^-{ a^+_C} & {\Sigma^{-1}\mathfrak{S}} \ar[u]_-{\Sigma^{-1}a^+_{SC}} \ar[l]^{\mathfrak{U}}}\] and
splice in the commutative square representing a morphism \( \Sigma ^{-1}f:\Sigma ^{-1}SC\rightarrow \Sigma ^{-1}D \)
(in \( \mathfrak{L}_{0}^{-1} \)) \[\xymatrix{{\Sigma^{-1}\coend((SC)^+)}\ar[r]^-{\Sigma^{-1}f}&{\Sigma^{-1}D^+}\\{\Sigma^{-1}\mathfrak{S}}\ar@{=}[r]\ar[u]^{\Sigma^{-1}a^+_{SC}}&{\Sigma^{-1}\mathfrak{S}}\ar[u]_{\Sigma^{-1}a_D}}\]The
``outer rim'' of the result is precisely \ref{dia:smallsuspensionsquare}.
\end{proof}
The following is straightforward, but should be said

\begin{defn}
Given a coalgebra \( C \) over an operad \( \mathfrak{R}_{1} \) with structure
map
\[
a:\mathfrak{R}_{1}\rightarrow \coend (C)\]
and a morphism of operads
\[
f:\mathfrak{R}_{2}\rightarrow \mathfrak{R}_{1}\]
 the \emph{pullback}, \( f^{*}C \) is defined to be the coalgebra over \( \mathfrak{R}_{2} \)
with structure map 
\[
a\circ f:\mathfrak{R}_{2}\rightarrow \coend (C)\]

\end{defn}
We have two distinct ways of converting objects of \( \mathfrak{L} \) into
objects of \( \Sigma ^{-k}\mathfrak{L} \): by desuspending, as in \ref{prop:desuspensioninj},
or by pulling back structure morphisms over 
\[
\mathfrak{U}^{k}=\Sigma ^{-k}\mathfrak{V}^{k}:\Sigma ^{-k}\mathfrak{S}\rightarrow \mathfrak{S}\]
as in \ref{cor:suspensiondown}. 

The first method desuspends the underlying chain-complex (so we get a chain-complex
that extends into negative dimensions), and the second does not. These two functors
are related by:

\begin{thm}
\label{th:kfoldsuspension}Let \( C \) and \( D \) be m-coalgebras over \( \mathfrak{S} \)
with structure-maps \( a_{C} \) and \( a_{D} \), respectively, let \( k \)
be an integer \( >0 \). Then \( S^{k}C \) is equivalent to \( S^{k}D \) in
\( \mathfrak{L} \) if and only if \( (\mathfrak{U}^{k})^{*}C^{+} \) is equivalent
to \( (\mathfrak{U}^{k})^{*}D^{+} \) in \( \Sigma ^{-k}\mathfrak{L} \).
\end{thm}
\begin{proof}
This follows by an inductive application of \ref{th:operadsuspensiondiagram}.
\end{proof}
\begin{defn}
\label{def:mathfrakminfty}Consider the inverse system of operads
\[
\cdots \rightarrow \Sigma ^{-i}\mathfrak{S}\rightarrow \Sigma ^{-i+1}\mathfrak{S}\rightarrow \cdots \rightarrow \mathfrak{S}\]
whose maps are 
\[
\Sigma ^{-i}\mathfrak{U}:\Sigma ^{-i-1}\mathfrak{S}\rightarrow \Sigma ^{-i}\mathfrak{S}\]
Denote the corresponding inverse limit by \( \mathfrak{S}_{-\infty } \). The
morphism \( \mathfrak{U}:\Sigma ^{-1}\mathfrak{S}\rightarrow \mathfrak{S} \)
induces an operad-endomorphism
\[
\mathfrak{U}:\mathfrak{S}_{-\infty }\rightarrow \mathfrak{S}_{-\infty }\]

In addition, \( \mathfrak{S}_{-\infty } \) comes equipped with canonical operad-morphisms
\[
\mathfrak{p}_{n}:\mathfrak{S}_{-\infty }\rightarrow \Sigma ^{-n}\mathfrak{S}\]

\end{defn}
\begin{rem}
We may think of \( \mathfrak{S}_{-\infty } \) as ``\( \Sigma ^{-\infty }\mathfrak{S} \)''
and the projections \( \mathfrak{p}_{n} \) as ``\( \mathfrak{U}^{\infty -n} \)''.
This makes some sense if one notes that all of the \( \mathfrak{U}^{n} \) maps
are \emph{surjective} so that there are isomorphisms
\[
\Sigma ^{-n}\mathfrak{S}\rightarrow F^{n}\{\Sigma ^{-i}\mathfrak{S}\}\]
where \( F^{n}\{\Sigma ^{-i}\mathfrak{S}\} \) is a ``finite'' portion of
the inverse limit --- the sub-operad of 
\[
\prod _{i=0}^{n}\Sigma ^{-i}\mathfrak{S}\subset \prod _{i=0}^{\infty }\Sigma ^{-i}\mathfrak{S}\]
composed of sequences \( \{\alpha _{0},\dots ,\alpha _{n}\} \) with \( \alpha _{i}=\mathfrak{U}(\alpha _{i+1}) \)
for all \( i<n \).
\end{rem}
\begin{defn}
\label{def:infinitedesuspL}Let \( \mathfrak{L}'_{-\infty } \) denote the category
of coalgebras over \( \mathfrak{S}_{-\infty } \) and let \( \mathfrak{L}_{-\infty }=\mathfrak{L}'_{-\infty }[H^{-1}] \),
where \( H \) is the class of morphisms inducing homology-isomorphisms. 

Let \( \Sigma ^{-\infty }\mathfrak{L}_{0}\subset \mathfrak{L}'_{-\infty } \)
be the subcategory with the property that each object, \( \mathrf {O} \), of
\( \Sigma ^{-\infty }\mathfrak{L}_{0} \) is of the form \( \mathfrak{p}_{k}^{*}Z \), 
\begin{enumerate}
\item where \( Z\in \Sigma ^{-k}\mathfrak{L}_{0} \) for some value of \( k \), and
\item \( \mathfrak{p}_{k}:\mathfrak{S}_{-\infty }\rightarrow \Sigma ^{-k}\mathfrak{S} \),
as defined in \ref{def:mathfrakminfty}.
\end{enumerate}
Given an object \( \mathrf {O} \), of \( \Sigma ^{-\infty }\mathfrak{L}_{0} \)
of the form \( \mathfrak{p}_{k}^{*}Z \), we will call the integer \( k \)
the \emph{level} of \( \mathrf {O} \). In addition, define \( \Sigma ^{-\infty }\mathfrak{L}=\Sigma ^{-\infty }\mathfrak{L}_{0}[H^{-1}] \),
where \( H \) is the class of morphisms inducing homology isomorphisms. Then
morphisms \( m:\mathrf {O}_{1}\rightarrow \mathrf {O}_{2} \) are diagrams of
the form 
\[
\mathrf {O}_{1}\leftrightarrow Z_{1}\leftrightarrow \cdots \leftrightarrow Z_{n}\leftrightarrow \mathrf {O}_{2}\]
with each \( \leftrightarrow  \) representing a morphism to the left or right,
and all left-pointing arrows inducing homology equivalences. We will define
the \emph{level of the morphism} \( m \) to be 
\[
\max \left\{ \mathrm{level}(\mathrf {O}_{1}),\mathrm{level}(Z_{1}),\dots ,\mathrm{level}(Z_{n}),\mathrm{level}(\mathrf {O}_{2})\right\} \]

\end{defn}
\begin{cor}
\label{cor:arbitrarysuspension}Let \( C \) and \( D \) be m-coalgebras. Then
there exists an integer \( k>0 \) and a morphism \( S^{k}C\rightarrow S^{k}D \)
in \( \mathfrak{L} \) if and only if there exists a morphism \( \mathfrak{p}_{0}^{*}C^{+}\rightarrow \mathfrak{p}_{0}^{*}D^{+} \)
in \( \Sigma ^{-\infty }\mathfrak{L} \). The two morphisms induce identical
maps of chain-complexes. In particular, there exists an integer \( k>0 \) such
that \( S^{k}C \) is equivalent to \( S^{k}D \) in \( \mathfrak{L} \) if
and only if \( \mathfrak{p}_{0}^{*}C^{+} \) is equivalent to \( \mathfrak{p}_{0}^{*}D^{+} \)
in \( \Sigma ^{-\infty }\mathfrak{L} \).
\end{cor}
\begin{proof}
We begin by proving the \emph{if}-part: If there exists a morphism \( S^{k}C\rightarrow S^{k}D \)
then, by \ref{th:kfoldsuspension}, \( (\mathfrak{U}^{k})^{*}C^{+}\cong (S^{k}C)^{+} \)
and \( (\mathfrak{U}^{k})^{*}D^{+}\cong (S^{k}D)^{+} \) in \( \Sigma ^{-k}\mathfrak{L} \).
Since \( \mathfrak{p}_{k}\circ \mathfrak{U}^{k}=\mathfrak{p}_{0} \), the morphism
\( S^{k}C\rightarrow S^{k}D \) induces a morhism 
\[
\mathfrak{p}_{0}^{*}C^{+}\rightarrow \mathfrak{p}_{0}^{*}D^{+}\]
in \( \mathfrak{L}_{-\infty } \).

One the other hand, suppose there exists a morphism \( \mathfrak{p}_{0}^{*}C^{+}\rightarrow \mathfrak{p}_{0}^{*}D^{+} \)
in \( \Sigma ^{-\infty }\mathfrak{L} \), of level \( k \) (see \ref{def:infinitedesuspL}).
Then there exists a sequence 
\begin{equation}
\label{eq:equivalencesequence1}
\mathfrak{p}_{0}^{*}C^{+}\leftrightarrow Z_{1}\leftrightarrow \cdots \leftrightarrow Z_{t}\leftrightarrow \mathfrak{p}_{0}^{*}D^{+}
\end{equation}
where each arrow, \( \leftrightarrow  \), points to the left or right and the
left-pointing arrows induce homology isomorphisms. Suppose \( Z_{i}=\mathfrak{p}_{k_{i}}^{*}W_{i} \),
where \( W_{i}\in \Sigma ^{-k_{i}}\mathfrak{L}_{0} \). Then \( k=\max \{k_{i}\} \).
We convert \ref{eq:equivalencesequence1} into a sequence
\[
(\mathfrak{U}^{k})^{*}C^{+}\leftrightarrow (\mathfrak{U}^{k-k_{1}})^{*}Z_{1}\leftrightarrow \cdots \leftrightarrow (\mathfrak{U}^{k-k_{t}})^{*}Z_{t}\leftrightarrow (\mathfrak{U}^{k})^{*}D^{+}\]
which defines a morphism \( (\mathfrak{U}^{k})^{*}C^{+}\rightarrow (\mathfrak{U}^{k})^{*}D^{+} \)
in \( \Sigma ^{-k}\mathfrak{L} \). The conclusion follows by \ref{th:kfoldsuspension}.
\end{proof}
\begin{cor}
\label{cor:arbitrarysuspensionspaces}Let \( X \) and \( Y \) be pointed,
simply-connected semisimplicial sets. Then there exists an integer \( k\geq 0 \)
and a pointed map of topological realizations \( f:|S^{k}X|\rightarrow |S^{k}Y| \)
if and only if there exists a morphism \( g:\mathfrak{p}_{0}^{*}\mathscr{C}(X)^{+}\rightarrow \mathfrak{p}_{0}^{*}\mathscr{C}(Y)^{+} \)
in \( \Sigma ^{-\infty }\mathfrak{L} \). In this case, the following diagram
commutes in \( \mathfrak{L} \) \begin{equation}\xymatrix{{\Sigma^k(\mathfrak{U}^k)^*\cfp{X}}\ar@{=}[d] \ar[r]^{\Sigma^kg}& {\Sigma^k(\mathfrak{U}^k)^*\cfp{Y}} \ar@{=}[d]\\ {\cfp{S^kX}} \ar[d]_{\iota_{S^kX}}& {\cfp{S^kY}}\ar[d]^{\iota_{S^kY}} \\ {\mathscr{C}(\ddelta(|S^kX|))}\ar[r]_{\mathscr{C}(\ddelta(f))} & {\mathscr{C}(\ddelta(|S^kY|))}}\label{dia:spacecorollary}\end{equation} where
\begin{enumerate}
\item \( \Sigma ^{k} \) denotes suspension of chain-complexes and base-operads. This
is suspension in the sense of dimension-shifting.
\item \( \ddelta (*) \) denotes the 2-reduced singular complex functor.
\item \( |*| \) denotes topological realization functor.
\item \( \iota _{X}:\mathscr{C}(X)\rightarrow \mathscr{C}(\ddelta (|X|)) \) is the
canonical (injective) equivalence.
\end{enumerate}
In particular \( X \) and \( Y \) are stably homotopy equivalent if and only
if \( \mathfrak{p}_{0}^{*}\cfp{{X}} \) and \( \mathfrak{p}_{0}^{*}\cfp{{Y}} \)
are equivalent in \( \Sigma ^{-\infty }\mathfrak{L} \).

\end{cor}
\begin{rem}
Note that \( \mathfrak{L} \) is the localized category \( \mathfrak{L}_{0}[H^{-1}] \),
so that the commutativity of \ref{dia:spacecorollary} in \( \mathfrak{L} \)
implies its \emph{homotopy} commutativity in \( \mathfrak{L}_{0} \). 

Note that the basepoints have fallen out of the picture: the stable homotopy
type of \( X \) is determined by the reduced m-coalgebra \( \cfp{{X}} \).
\end{rem}
\begin{proof}
This follows from \ref{cor:arbitrarysuspension}, which implies that \( S^{k}\cf{{X}} \)
is equivalent to \( S^{k}\cf{{Y}} \) in \( \mathfrak{L} \) and \cite[Theorem 9.9]{Smith:2000}.
\end{proof}

\providecommand{\bysame}{\leavevmode\hbox to3em{\hrulefill}\thinspace}

\end{document}